\documentclass{article}
\usepackage{amsthm}
\usepackage{amsfonts}
\usepackage{cite}
\usepackage{amsmath, amscd}
\newtheorem{theorem}{Theorem}[section]
\newtheorem{definition}[theorem]{Definition}

\newtheorem{lemma}[theorem]{Lemma}
\newtheorem{example}[theorem]{Example}

\begin{document}
\author{Jeremy Berquist}
\title{Singularities on Demi-Normal Varieties}
\maketitle

\noindent
\linebreak
\textbf{Abstract.  }   The birational classification of varieties inevitably leads to the study of singularities.  The types of singularities that occur in this context have been studied by Mori, Koll\'ar, Reid, and others, beginning in the 1980s with the introduction of the Minimal Model Program.   Normal singularities that are terminal, canonical, log terminal, and log canonical, and their non-normal counterparts, are typically studied by using a resolution of singularities (or a semi-resolution), and finding numerical conditions that relate the canonical class of the variety to that of its resolution.  In order to do this, it has been assumed that a variety $X$ is has a $\mathbb{Q}$-Cartier canonical class:  some multiple $mK_X$ of the canonical class is Cartier.  In particular, this divisor can be pulled back under a resolution $f: Y \rightarrow X$ by pulling back its local sections.  Then one has a relation $K_Y \sim \frac{1}{m}f^*(mK_X) + \sum a_iE_i$.  It is then the coefficients of the exceptional divisors $E_i$ that determine the type of singularities that belong to $X$.  It might be asked whether this $\mathbb{Q}$-Cartier hypothesis is necessary in studying singularities in birational classification.  In \cite{dFH09}, de Fernex and Hacon construct a boundary divisor $\Delta$ for arbitrary normal varieties, the resulting divisor $K_X + \Delta$ being $\mathbb{Q}$-Cartier even though $K_X$ itself is not.  This they call (for reasons that will be made clear) an $m$-compatible boundary for $X$, and they proceed to show that the singularities defined in terms of the pair $(X,\Delta)$ are none other than the singularities just described, when $K_X$ happens to be $\mathbb{Q}$-Cartier.  Thus, a wider context exists within which one can study singularities of the above types.  In the present paper, we extend the results of \cite{dFH09} still further, to include demi-normal varieties without a $\mathbb{Q}$-Cartier canonical class.  Our secondary aim is to give an introduction to the study of those non-normal varieties that appear naturally in the Minimal Model Program, culminating with the extension of de Fernex and Hacon's results to these varieties.

\newpage
\tableofcontents
\addcontentsline{}{}{}
\newpage
\begin{section}{Introduction}

Let us first recall briefly the birational classification of surfaces.  Birationally, a surface $X$ is nonsingular and projective (by theorems of Nagata, Chow, and Hironaka).  Certain rational curves on $X$ (the (-1)-curves) can be contracted, giving a birational morphism $X \rightarrow X'$ to another (nonsingular and projective) surface.  This process stops after finitely many steps, yielding a ``relatively minimal model" of $X$.  Dealing with the rational and ruled surfaces as separate cases, there exists a unique such model, called the minimal model of $X$.  It is characterized uniquely by the property that $K_{X'}$ is nef (numerically effective).

The Minimal Model Program (MMP, or Mori's program) extends the birational classification of surfaces to higher dimensions; that is, it looks for a way to find a birational model $X'$ where $K_{X'}$ is nef.  It is complete in dimension three, and large parts are complete in all dimensions.  It is hoped that having a nef canonical divisor will allow further useful steps in birational classification.

In fact, on a different scale, MMP (if true) is used in moduli theory to reduce to cases where the canonical class is generated by its global sections.   Then one reduces to cases of ``general type," where there exists a birational morphism from $X$ to its canonical model $Y := \textnormal{Proj }\oplus_{m \geq 0} H^0(X, mK_X).$  Then the Hilbert polynomial of the canonical divisor of $Y$ is the principle discrete invariant used in constructing moduli spaces.

Birationally, $X$ is smooth and projective.  However, MMP introduces singularities as it finds a nef birational model.  Thus the canonical model may have some singularities.  These singularities are well-understood in certain cases.  We attempt to resolve issues dealing with the remaining cases.  

Letting $Y$ be the canonical model of $X$, $K_Y$ is used to measure how singular the canonical model is.  From a resolution of singularities $f: Y' \rightarrow Y$, one looks at the coefficients of the exceptional divisors in $K_{Y'} \sim f^*K_Y + \Sigma a_iE_i.$  The discrepancy is the minimum of the coefficients $a_i$.  Depending on whether the discrepancy is at least -1, greater than -1, at least 0, or greater than 0, the canonical model has log canonical, log terminal, canonical, or terminal singularities.

Two issues immediate arise.  First, if we want to compare divisors on $Y$ with those on $Y'$, we need $f$ to be an isomorphism over the codimension one points of $Y$.  For varieties where $K_Y$ makes sense (in particular, assuming Serre's $S_2$ property), this is the same thing as saying that $Y$ is  normal.  Second, providing that we can talk reasonably about $K_Y$, we still need to be able to pull it back.  This is possible if $K_Y$ is $\mathbb{Q}$-Cartier, and until recently this has been a standard assumption.  

These two assumptions can be relaxed independently of one another in order to deal with more general varieties.  That is, one might want to deal with non-normal varieties, and one might want to dispense with the $\mathbb{Q}$-Cartier hypothesis on $K_Y$.   The non-normal varieties studied in the context of MMP are the demi-normal varieties.  These are varieties with properties $G_1, S_2,$ and seminormality.  (Note that a normal variety is one that is $R_1$ and $S_2$, and that $G_1$, or Gorenstein in codimension one, is slightly weaker than $R_1$.  We will discuss these properties in the second section of this paper.)  Equivalently, these are varieties with Serre's $S_2$ property, and such that the codimension one singularities are of the simplest possible type, being double normal crossings.  For the second assumption, de Fernex and Hacon have solved the problem of pulling back $K_Y$ on a normal variety when $K_Y$ is not $\mathbb{Q}$-Cartier.  In fact, they construct a boundary divisor $\Delta$ called an $m$-compatible boundary and study the singularities of the pair $(Y, \Delta)$; they also show that the resulting class of singularities corresponds with the usual class when $K_Y$ is $\mathbb{Q}$-Cartier.

It is our primary goal in this paper to resolve this two issues simulataneously.  In other words, we follow the lead of de Fernex and Hacon and dispense with a  $\mathbb{Q}$-Cartier canonical class, in our case treating demi-normal varieties.

There are technical reasons why we cannot deal with arbitrary non-normal varieties.  For instance, to be able to talk about Weil divisors on a non-normal variety, we need to assume properties $G_1$ and $S_2$.  Furthermore, if we plan on somehow resolving the singularities of a non-normal variety, then it seems we need to assume seminormality as well.   Hironaka's resolution of singularities is too strong for non-normal varieties, at least if we want an isomorphism over the codimension one structure.  Koll\'ar solves this problem by showing that a semiresolution exists for demi-normal varieties.  There are some singularities, but they are only of two types and easily described.  In fact, the only singularities on a semiresolution are double points, either double normal crossings (locally, solutions of $xy = 0$) or pinch points (local solutions of $x^2 - y^2z = 0$).

We resolve both issues by following de Fernex and Hacon with Koll\'ar's semiresolutions in mind.  The main result in this direction is the following generalization of a theorem of de Fernex and Hacon:

\begin{theorem}  Let $X$ be a demi-normal variety.  Then for any semiresolution $f: Y \rightarrow X$, $m$-compatible boundaries exist with respect to $f$, for all sufficiently divisible positive integers $m$.
\end{theorem}

We first define a generalized pullback that works for alll Weil divisors on a demi-normal variety.  Then we show that the resulting theory is equivalent to the theory of pairs.  The crux of the matter is that we can assume that a boundary divisor exists for $X$, with the property that the singularities of the resulting pair (as they are already understood) are the same as the generalized singularities of de Fernex and Hacon, interpreted for demi-normal varieties.

There is some surprising behavior in the resulting classes of singularities.  In fact, we show that following definitions properly leads to the existence of a semi-canonical variety that is not sem-log terminal.  This is evidently impossible when dealing with varieties with a $\mathbb{Q}$-Cartier canonical class.  Nevertheless, our example shows how new characteristics arise when relaxing this hypothesis.

The paper is organized as follows.  In the seond section, we review semiresolutions of demi-normal varieties and discuss why ``demi-normal" is the correct notion of a non-normal variety in the context of MMP.  Next, we present the theory of divisors on varieties with properties $G_1$ and $S_2$.  It turns out that if we consider only rank one reflexive sheaves, which correspond to Weil divisors on this class of varieties, we obtain enough flexibility to pull back arbitrary (not necessarily $\mathbb{Q}$-Cartier) Weil divisors.  Working with sheaves turns out to be easier than working with divisors.  We also discuss the standard formation of the canonical or dualizing class on a variety with the $G_1$ and $S_2$ properties.  This is essential to developing a definition of singularities that extends the standard definitions in MMP.  In the fourth section, we give our proof of Theorem 1.1.  In the following section, we give new definitions of singularities of demi-normal varieties without a $\mathbb{Q}$-Cartier canonical class, following the definitions of de Fernex and Hacon for normal varieties without such a class.  Finally, we give an example of a variety that is semi-canonical but not semi-log terminal.  It is to be wondered if a different definition of singularities is perhaps better, or whether this type of unexpected behavior is a result of properties of arbitrary demi-normal varieties themselves.
\end{section}

\begin{section}{Semiresolutions}
This section contains definitions of properties assumed to hold on the varieties that we study.  At the end, we work out an example of a semiresolution, demonstrating the standard technique of finding a semiresolution by finding a resolution of the normalization and gluing along the birational transform of the conductor.

We assume that $X$ is a reduced, equidimensional, quasi-projective variety over an algebraically closed field of characteristic zero.  In particular, the normalization is finite and resolution of singularities is true for $X$.

We also assume three conditions on the local rings of $X$:  (i) $S_2$, (ii) $G_1$, and (iii) $SN$.  Equivalently, $X$ is $S_2$ and there is an open subvariety $U$, whose complement has codimension at least two, such that for any closed singular point $x \in U$, $\mathcal{O}_{X,x}^* \cong k[[x_1, ...,x_n]]/(x_1x_2),$ where the star denotes completion.  Such varieties are called demi-normal.

 Next, we will define these properties and give examples to show how they are distinct and related to one another.  

We say that a local ring is $S_n$ if its depth at the maximal ideal is greater than the minimum of $n$ and its dimension:  $$\textnormal{depth } \mathcal{O}_{X,x} \geq \textnormal{min } \{n, \textnormal{dim } \mathcal{O}_{X,x} \}.$$   It is Cohen-Macaulay if it is $S_n$ for all $n$, or equivalently (noting that the depth is always bounded above by the dimension), if the depth is equal to the dimension.  In more advanced terms, the dualizing complex exists for our varieties, and being Cohen-Macaulay means that the dualizing complex is quasi-isomorphic to a complex concentrated in a single term.  We'll say more about the dualizing complex and dualizing sheaf in the following sections.  The condition $S_2$ in particular has the geometric interpretation that a function regular in codimension one is regular everywhere.  It is a base condition that allows us to compare $X$ with its normalization by looking just at the codimension one structure of $X$.    

A local ring is Gorenstein if it is Cohen-Macaulay and the dualizing sheaf (meaning in this case the single term in the dualizing complex) is invertible.  Any smooth point has a corresponding local ring that is Gorenstein.  In fact, we have the following implications:  $$\textnormal{complete intersection } \implies \textnormal{Gorenstein } \implies \textnormal{Cohen-Macaulay }.$$  So hypersurfaces are also Gorenstein, and we have a large class of singular, Gorenstein local rings given by the local rings of hypersurface singularities.  Also, any Gorenstein variety (one whose local rings are all Gorenstein) trivially has $S_2$.  The condition $G_1$ means Gorenstein in codimension one; in other words, at any codimension one or codimension zero point, the local ring is Gorenstein.

We say that a ring extension $A \hookrightarrow B$ is subintegral if it is finite, a bijection on prime spectra, and each residue field extension $k(p) \hookrightarrow k(q)$ is an isomorphism (here $k$ is the base field described above).  For any finite extension $A \hookrightarrow B$, there exists a maximal subintegral extension of $A$ in $B$, called the seminormalization of $A$ in $B$.  We say that $A$ is seminormal (SN) in $B$ if it equals its seminormalization in $B$.  We call a reduced ring $A$ seminormal if it is seminormal in its normalization.  A variety $X$ is seminormal if its local rings are all seminormal.  Equivalently, any proper bijection to $X$ is an isomorphism.

Of course, any normal variety is $S_2$, $G_1$, and seminormal.   In fact, normality is characterized by the two conditions $S_2$ and $R_1$, or regular in codimension one, and clearly $R_1$ is stronger than $G_1$ and every normal ring is seminormal.  The typical model for a demi-normal variety is that of a simple normal crossing divisor in a smooth variety.  We think of demi-normal varieties as being of the very simplest singular type in codimension one, with regular functions in codimension one being regular everywhere.  Several other examples will be given throughout this paper.

The two conditions $S_2 + G_1$ together imply that rank one reflexive sheaves determine subschemes of pure codimension one in $X$ (Weil divisors in the normal case).  The conditions $S_2 + SN$ imply that the normal locus is a reduced subscheme of pure codimension one.  Finally, the condtions $G_1 + SN$ imply that a each codimension one localization $R$, the maximal ideal $m$ is the radical of the normalization $\overline{R}$, and the $R/m$-dimension of $\overline{R}/m$ is exactly two.

\begin{example}  For a variety that is $S_2$ and $G_1$ but not seminormal, take any hypersurface that is not seminormal.  For instance, the cusp $k[x,y]/(x^2 - y^3)$.  Its normalization is the affine line, and maps finitely and bijectively onto the cusp, but is clearly not an isomorphism. 
\end{example}

\begin{example}   For a variety that is $S_2$ and SN but not $G_1$, consider the three coordinate axes in 3-dimensional affine space.  This is reduced and one-dimensional, hence $S_2$ (just $S_1$ is sufficient), and seminormal by \cite{LV81}.  However, there are not double normal crossing singularities only, so it cannot be Gorenstein.  
\end{example}

\begin{example}  Finally, for a variety that is $G_1$ and SN but not $S_2$, consider two planes in 4-space meeting at a point.  The coordinate ring is $k[x,y,z,w]/((x,y) \cap (z,w)).$  This is $G_1$ because it is smooth except at the origin, which is a point of codimension two.  SN can be proved directly using the following condition:  a ring $A$ is seminormal in $B$ if for every $b \in B$, $b^2, b^3 \in A$ implies $b \in A$.  It is not $S_2$ because at the origin, the element $x+z$ forms a maximal regular sequence:  the depth at the origin is exactly one.
\end{example}

We call a variety $X$ semismooth if for any closed singular point $x$, the completion $\mathcal{O}_{X,x}^*$ is isomorphic to either $k[[x_1, \ldots, x_n]]/(xy)$ or $k[[x_1, \ldots, x_n]]/(x_1^2-x_2^2x_3).$  The latter type of singularity is called a pinch point.  Pinch points are notorious because blowing up at the origin creates another pinch point; the only way to resolve its singularities is to blow up the entire double locus.  Locally away from the origin, the local ring of a pinch point has double normal crossing singularities, and a pinch point in dimension two can be obtained from a double normal crossing point as a quotient singularity.  Both these singularities are simple examples of partition singularities.  See \cite{vS87}.

The normalization of a semismooth variety is smooth.  The double locus of $X$ and its preimage in the normalization $\overline{X}$ are both smooth, and the induced morphism on the double loci is a double cover, ramified along the pinch locus.  Note that a variety is smooth if all its closed points are smooth (that is, the local rings are all regular local rings), because the localization of a regular local ring is regular.  Also, a local ring is regular if and only if its completion is also a regular local ring.  This is why we can restrict our attention to completions of local rings at closed points in the definition of semi-smoothness. 

Given a demi-normal variety $X$ and an open subvariety $U$ where the only singularities are double normal crossings, there is a morphism $f: Y \rightarrow X$ with the following properties:  
\begin{itemize}
\item{$f$ is proper and birational},
\item{$f|_{f^{-1} U}$ is an isomorphism},
\item{$Y$ is semismooth},
\item{no component of the double locus of $Y$ is exceptional}.
\end{itemize}

We call such a morphism a semiresolution of $X$.  The existence of $f$ (in fact, a stronger form of semiresolution, where we incorporate divisors on $X$, also exists, given reasonable conditions on the divisor) has been proved by Koll\'ar in his paper ``Semi Log Resolutions."  Here $Y$ is not normal, but every exceptional divisor determines a discrete valuation ring at its generic point.  This enables us to pull back Weil divisors on $X$ (subschemes determined by rank one reflexive sheaves) via semiresolutions by considering the order of vanishing in each such discrete valuation ring.

Since $Y$ is not smooth, there is not everywhere a local system of parameters.  However, $\overline{Y}$ is smooth, and we obtain our model for a normal crossing divisor on $Y$ as the gluing along the double locus of a normal crossing divisor on $\overline{Y}$ that makes transverse intersections with the double locus.  Equivalently, the local models for a global normal crossing divisor $D$ on a semismooth variety $Y$ are as follows:

When $Y$ is given by $x_1x_2=0$ in $\mathbb{A}^n_{x_1, \ldots, x_n}$, the local model for $D$ is $D = \Pi_{i \in I} x_i = 0$, for some set $I \subseteq \{3, \ldots, n \}$.  For the pinch point, given locally (analytically) by $x_1^2-x_2^2x_3 = 0$, the local model is $D = (\Pi_{i \in I} x_i = 0) + D_2$, for some $I \subseteq \{4, \ldots, n \}$ and where either $D_2 = 0$ or $D_2 = (x_1 = x_3 = 0).$  

The analog of a log resolution of a pair is a semiresolution with the property that the exceptional divisors, the double locus, and the birational transform of the given divisor form a global normal crossing divisor.  When $X$ is demi-normal with an open subvariety $U$ as described above, and $D$ is a divisor such that $D|_U$ is smooth and disjoint from the singular locus of $U$, a semi log resolution of $(X,D)$ exists.  The idea is to pass to the normalization, find a log resolution of singularities, and then glue along the birational transform of the double locus.

The main tool in constructing semiresolutions is the universal pushout.  Given a variety $Y$, a closed subscheme $B$, and a finite morphism $B \rightarrow B/ \tau$ (we use the language of a quotient singularity, because that is often the form of a given finite morphism from a closed subscheme of $Y$), there is a commutative diagram $$\begin{CD}
B           @>>>     Y \\
@VVV                    @VVV \\
B/\tau   @>>>     X \\
\end{CD}$$
where both horizontal arrows are inclusions of closed subschemes.  The resulting morphism $Y \rightarrow X$ to the universal pushout is finite.  It is also birational; in fact, it agrees with $B \rightarrow B/\tau$ on $B$ and is an isomorphism elsewhere.  The pushout occurs naturally in connection with the normalization of a variety.  In fact, when $X$ is a variety and $Y$ is its normalization, then one obtains $X$ from $Y$ as a universal pushout by gluing along the conductor loci.  (When we say gluing we just mean constructing the pushout; $Y$ is ``glued" to $X$ along the morphism $B \rightarrow B/\tau.$)  The universal pushout does indeed have a universal property.  In fact, given morphisms $B/\tau \rightarrow X_0$ and $Y \rightarrow X_0$ that pull back to the same morphism from $B$, there exists a unique morphism $X \rightarrow X_0$ commuting with the given morphisms.  Finally, the universal pushout and pullback are related in the usual scheme-theoretic way:  given affine schemes, their pushout is the affine scheme whose coordinate ring is the pullback of the coordinate rings of the given schemes.

We close this section with an example of a semiresolution.

\begin{example}  Consider the surface defined by $k[x,y,z]/(xy)$, and let $\mathbb{Z}_2$ act on it by $x \mapsto -x, y \mapsto -y, z \mapsto -z$.  The resulting quotient is a non-normal variety $X$ with two normal components.  A semiresolution is obtained by resolving the components and gluing along the birational transform of the line of intersection.
\end{example}  
By definition, the hyperquotient singularity $X$ is obtained as a residue class of a ring of invariants.  If $\mathbb{Z}_2$ acts on $k[x,y,z]$ as stated, the ring of invariants is $k[x^2,xy,y^2,xz,yz,z^2]$.  We get the coordinate ring of $X$ by annihilating the intersection of the ideal $(xy)$ with this ring.  The intersection ideal has generators $xy,x^2y^2,xyz^2,x^2yz,$ and $xy^2z$.  Therefore, in terms of generators and relations, the coordinate ring $\mathcal{O}_X$ is given by $$k[u_0,u_1,u_2,u_3,u_4,u_5]/(I + J),$$ where $$I = (u_0u_2-u_1^2, u_2u_5-u_4^2, u_0u_5-u_3^2,u_1u_3-u_0u_4, u_3u_4-u_1u_5)$$ and $$J = (u_1,u_0u_2,u_3u_4,u_0u_4,u_2u_3).$$  If we simplify, then we obtain $k[u_0,u_2,u_3,u_4,u_5]$ modulo the ideal $$(u_2u_5-u_4^2,u_0u_5-u_3^2,u_0u_2,u_3u_4,u_0u_4, u_2u_3).$$
Note that the spectrum of this ring has two components, given by $(u_0=u_3=0)$ and $(u_2=u_4=0)$.  Each of these defines a quadric cone, and the cones are identified along the line $\mathbb{A}^1 \cong \textnormal{Spec } k[u_5]$.  We can obtain the same ring as a pullback, where the maps $c$ and $d$ are the quotient maps:
$$\begin{CD}
\mathcal{O}_X                @>>>   k[u_0,u_3,u_5]/(u_0u_5-u_3^2)\\
@VVV                                @VcVV \\
k[u_2,u_4,u_5](u_2u_5-u_4^2) @>d>>  k[u_5]
\end{CD}$$  
If we compute the $\mathbb{Z}_2$-quotient of each component of $k[x,y,z]/(xy)$ and then glue, we get $\mathcal{O}_X$.  In fact, the quadric cone is a $\mathbb{Z}_2$-quotient of affine two-space $\mathbb{A}^2$.

The components of $X$ are both normal, and each of their singularities is resolved by blowing up the origin.  In each of these two blowups, the birational transform of the intersection line $\mathbb{A}^1$ appears in only one chart.  Gluing over the corresponding ring $k[u_5']$, we get the pullback diagram 
$$\begin{CD}
\mathcal{O}_Y @>>>    k[u_3',u_5'] \\
@VVV                  @Vc'VV \\
k[u_4',u_5']  @>d'>>  k[u_5']
\end{CD}$$
Here $\mathcal{O}_Y$ is $k[u_3',u_4',u_5']/(u_3'u_4')$.  On the other charts, no gluing takes place and the variety remains smooth.  Thus $Y$ has only double normal crossing singularities and is semismooth.  

Globally, $Y$ is the pushout 
$$\begin{CD}
\tilde{\mathbb{A}}^1 @>>> Y_2 \\
@VVV                      @VVV \\
Y_1                  @>>> Y
\end{CD}$$
Here $Y_1$ and $Y_2$ are the resolutions of the components $X_1 := (u_0=u_3=0)$ and $X_2 := (u_2 = u_4 =0)$ of $X$, and $\tilde{\mathbb{A}}^1$ is the birational transform of the intersection of $X_1$ and $X_2$.  The universal property of the pushout implies that $f: Y \rightarrow X$ exists.  The open sets of $X$ can be identified with pairs of open sets that have the same pullback to $\mathbb{A}^1$.  Then $Y \rightarrow X$ is an isomorphism over $(U_1,U_2)$, where $U_i$ is the complement of the origin in $X_i$.   By construction, the singular locus of $Y$ is the birational transform $\mathbb{A}^1$, which is not exceptional.

It needs to be checked that $f$ is proper.  This can be done using the valuative criterion of properness.
\end{section}
\begin{section}{Pullbacks of Divisors}
Most of the discussion in this section follows \cite{Hart94} and \cite{dFH09}.  The only real novelty comes in combining those results to develop a theory of pullbacks of divisors on a demi-normal variety under their semiresolutions.

The two conditions $S_2$ and $G_1$ imply that rank one reflexive sheaves correspond to subschemes of pure codimension one.  If $\mathcal{I}$ is a rank one reflexive ideal sheaf, there is an exact sequence $$0 \rightarrow \mathcal{I} \rightarrow \mathcal{O}_X \rightarrow \mathcal{O}_Y \rightarrow 0,$$ where $Y$ is the subscheme determined by $\mathcal{I}$.  The point is that on a scheme with $S_2$ and $G_1$, a sheaf is reflexive if and only if it is $S_2$ (the same definition as for rings, except the depth of a module over a local ring is considered); thus forming the reflexive hull is an $S_2$-ification.  Then the depth at each local ring of dimension at least two is at least two.  The long exact sequence in local cohomology begins with $$0 \rightarrow H^0_x(X, \mathcal{I}_x) \rightarrow H^0_x(X, \mathcal{O}_{X,x}) \rightarrow H^0_x(X, \mathcal{O}_{Y,x}) \rightarrow H^1_x(X, \mathcal{I}_x) \rightarrow \cdots.$$  Then the vanishing of the local cohomology modules in dimensions smaller than the depth implies that $Y$ has associated points of pure codimension one.  In other words, $Y$ is a divisor.  An arbitrary rank one reflexive sheaf embeds in the sheaf of total quotient rings, and can be multiplied by an invertible sheaf so that it becomes an ideal sheaf.  Thus rank one reflexive sheaves correspond to divisors, where we interpret some of the coefficients as negative.

If $f: Y \rightarrow X$ is a birational morphism, and $D$ is a divisor on $X$, then the reflexive hull of the inverse image ideal sheaf, or $(\mathcal{O}_X(-D)\cdot \mathcal{O}_Y)^{\vee\vee}$ is again of rank one and reflexive, and the natural pullback, denoted $f^{\flat}D$, is determined by setting $\mathcal{O}_Y(-f^{\flat}D)$ equal to this sheaf.  The natural pullback is additive if one of the summands is Cartier.  We let addition of divisors correspond to the reflexive hull of the product of the sheaves, the product being taken in the sheaf of total quotient rings.  However, we don't have additivity when one of the summands is only $\mathbb{Q}$-Cartier.  This is the same problem that occurs when we try to pull back Weil divisors on normal varieties.  The group laws of the group of divisors are not respected, and so we need to come up with a more general type of pullback that preserves the group structure.  To remedy this problem, we make use of discrete valuation rings corresponding to generic points of divisors on a normal variety.  When $X$ is only $G_1$, not every divisor corresponds to a discrete valuation ring; however, most of our birational morphisms are semiresolutions, and so every exceptional divisor does in fact correspond to such a ring (recall that no singular component on a semiresolution is ever exceptional).

Let $v$ be a discrete valuation ring associated to a codimension one point of $X$ (if such a point exists), and $\mathcal{J}$ and coherent $\mathcal{O}_X$-submodule of $\mathcal{K}_X$, where $\mathcal{J}$ is nonzero on the unique component of $X$ containing the center of $v$.  The valuation $v(\mathcal{J})$ of $\mathcal{J}$ along $v$ is defined to be $$v(\mathcal{J}) := \textnormal{min} \{v(\phi) : \phi \in \mathcal{J}(U \cap X_i), U \cap c_{X_i}(v) \neq 0 \}.$$  In particular, when $v$ is associated to an exceptional divisor $E$ of a semiresolution, the valuation of an invertible sheaf is the coefficient of $E$ in the corresponding Cartier divisor.

Let $v$ be a valuation of a component $X_i$ of $X$, and let $m$ be any positive integer.  Let us temporarily write $v^{\flat}(D)$ for the valuation along a divisor, described above for the rank one reflexive sheaf corresponding to $D$.   If $\mathcal{J}$ is any rank one reflexive sheaf, associated to a divisor $D$, and nonzero on $X_i$, then $$\inf_{k \geq 1} \frac{v^{\flat}(kD)}{k} = \lim \inf_{k \rightarrow \infty} \frac{v^{\flat}(kD)}{k} = \lim_{k \rightarrow \infty} \frac{v^{\flat}(k!D)}{k!} \in \mathbb{R}.$$  The reason the limit is finite is that $\mathcal{O}_X(D) \subseteq \mathcal{O}_X(C)$ for some Cartier divisor $C$ when $X$ is quasi-projective.  This suggests a definition, as in \cite{dFH09}:  for any divisor $D$ and a discrete valuation ring associated to an appropriate codimension one point, the valuation of $v$ along $D$ is $$v(D) := \lim_{k \rightarrow \infty} \frac{v^{\flat}(k!D)}{k!} \in \mathbb{R}.$$  

This allows us to define a pullback under a semiresolution $f: Y \rightarrow X$:  the pullback of a divisor $D$ on $X$ is the divisor $$f^*D := f^{-1}_*D + \Sigma v_E(D) \cdot E,$$ where $f^{-1}_*D$ is the birational transform of $D$ and the sum is over the exceptional divisors $E$, with the valuation $v_E$ corresponding to $E$ as just described.  

This generalized pullback is additive if one of the summands is $\mathbb{Q}$-Cartier.  In particular, it generalizes the notion of the pullback of a $\mathbb{Q}$-Cartier divisor; in that case, we chose a multiple $mD$ that is Cartier, pulled back the local equations, and then divided by $m$.  We will be able to use this generalized pullback to define singularities on demi-normal varieties without a $\mathbb{Q}$-Cartier canonical class.  Again, it should be noted that we measure the singularity of $X$ by comparing the canonical classes of $X$ and that of one of its semiresolutions.  To make this comparison, it is imperative that $K_X$ can be pulled back under a semiresolution, and for our definitions to agree with the usual case where $K_X$ is $\mathbb{Q}$-Cartier, we need the pullback to be the same as the usual pullback when the canonical class is $\mathbb{Q}$-Cartier.  

This prompts a question:  what exactly is the canonical class $K_X$?  The most general definition is as follows:  a quasi-projective variety over $k$ possesses a dualizing complex.  Then the lowest cohomology of this complex is a coherent sheaf.  When $X$ is $S_2$, the so-called dualizing sheaf is also $S_2$  (this is a general fact that can be proved using the definition alone.)   If $X$ is also $G_1$, then by definition the dualizing sheaf is invertible in codimension one, and it is reflexive since on varieties with $S_2$ and $G_1$, a sheaf is $S_2$ if and only if it is reflexive.  Thus the dualizing sheaf corresponds to a divisor (actually, a class of divisors, equivalence classes being defined up to multiples by elements of $\mathcal{K}_X$).  We call the correponding divisor class $K_X$.  A more practical definition is to form the dualizing sheaf on an open subvariety $U$ containing all the codimension one points, that being defined by differential forms for hypersurface singularities, push the sheaf forward to $X$, and take the corresponding divisor class on $X$.  In fact, with $S_2$ and $G_1$, a reflexive sheaf is completely determined by its behavior in codimension one.  The point is that we know how to form the canonical class on any smooth variety, or any variety with only hypersurface singularities, such as the double normal crossing locus in a demi-normal variety, and that with appropriate conditions on $X$, we can push forward to all of $X$

We now have all the ingredients we need to state and prove Theorem 1.1.  We do this in the next section, and in the following section, we make our definitions of singularities on demi-normal varieties without a $\mathbb{Q}$-Cartier canonical class.  There are three relative canonical divisors that will be useful to us; here $Y \rightarrow X$ is a semiresolution of $X$:

\begin{itemize}
\item{the $m$-th limiting relative canonical divisor $K_{m, Y/X} := K_Y - \frac{1}{m}f^{\flat}(mK_X)$,}
\item{the relative canonical divisor $K_{Y/X} := K_Y + f^*(-K_X)$,}
\item{for a boundary divisor $\Delta$ such that $K_X + \Delta$ is $\mathbb{Q}$-Cartier, and $\Delta_Y$ its birational transform, the log relative canonical divisor $$K^{\Delta}_{Y/X} := K_Y + \Delta_Y - f^*(K_X + \Delta).$$}
\end{itemize}

We may take the lim sup of the coefficients of $K_{m,Y/X}$ to obtain the divisor $K^{-}_{Y/X}.$  Then for each $m$ such that $m(K_X + \Delta)$ is Cartier, and all sufficiently large $q$, we have $$K^{\Delta}_{Y/X} \leq K_{m,Y/X} \leq K_{mq, Y/X} \leq K^{-}_{Y/X} \leq K_{Y/X}.$$  See \cite{dFH09} for a proof of these facts in the normal case; the demi-normal case is proved in the same way.
\end{section}

\begin{section}{Existence of $m$-Compatible Boundaries}
In this section, we prove Theorem 1.1.  In order to do so, we must first show that semi log resolutions exist for suitable pairs.  Koll\'ar discusses existence of semi log resolutions in his paper ``Semi Log Resolutions."  In particular, given a demi-normal variety $X$ with double normal crossing locus $U$, a divisor $D$ must be smooth and disjoint from the singular locus on $U$ in order for semi log resolutions to exist.  In our case, we begin only with a collection of rank one reflexive sheaves, invertible in codimension one, and it is not immediately clear that we can obtain a semi log resolution for the associated divisor on $X$.  However, by choosing sections appropriately, we can arrange for a semi log resolution to exist.  In proving the main theorem, we follow the proof of \cite{dFH09}, modified suitably for demi-normal varieties.  The difficult part is to pass to the normalization, apply typical results for normal varieties, and then glue along the conductor.  It is important that the normalization of a semismooth variety is smooth.  In fact, given $X$ and a semiresolution $Y \rightarrow X$, the induced morphism on normalizations $\overline{Y} \rightarrow \overline{X}$ is a resolution of singularities.  This gives us the freedom to work with normal and smooth varieties and glue in order to obtain results on demi-normal and semismooth varieties.

We start with a lemma.

\begin{lemma}  Suppose that $\mathcal{J}$ is a divisorial sheaf, invertible in codimension one.  Then there exists a semi log resolution of $(X, \mathcal{J})$.
\proof  The point here is that the associated divisor needs to have a nice restriction to an open semismooth subscheme.  So we need to be able to choose a sufficiently general global section to make everything work.  

Choose a very ample sheaf $\mathcal{L}$ such that $\mathcal{L} \otimes \mathcal{J}^{\vee}$ is generated by global sections.  Dualizing via a nondegenerate global section shows that $\mathcal{L}^{\vee} \otimes \mathcal{J}$ is isomorphic to an ideal sheaf.  We start by blowing up this ideal sheaf.  Then $(\mathcal{L}^{\vee} \otimes \mathcal{J}) \cdot \mathcal{O}_Y$ is invertible.  In particular, since $\mathcal{L}$ is invertible, so is $\mathcal{J} \cdot \mathcal{O}_Y$.  The blowup is an isomorphism in codimension one by the hypothesis on $\mathcal{J}$.  Moreover, any subsequent morphisms will preserve the invertibility condition.  We would like to replace $\mathcal{J}$ by $\mathcal{J} \cdot \mathcal{O}_Y$, thereby assuming that $\mathcal{J}$ is invertible.  Let $f: Y \rightarrow X$ be this first blowup.  If $V \subset Y$ is an open set with complement of codimension at least two, then $f(V^c)$ is also closed and with codimension at least two.  Thus, we can shrink $U$ if necessary so that it contains the codimension one points and so that $f^{-1}(U) \subset V$.  Then if we obtain an isomorphism in codimension one starting from $Y$, then we will have the desired morphism by composing with $f$.

We denote by $U$ an open set in $X$ that is semismooth and contains the codimension one points of $X$, and we let $U^{sm}$ be the smooth locus of $U$.  Choose a very ample sheaf $\mathcal{L}$ such that $\mathcal{L} \otimes \mathcal{J}^{\vee}$ is globally generated.  Then $(\mathcal{L} \otimes \mathcal{J}^{\vee})|_{U^{sm}}$ is generated by images of global sections in $\Gamma(X, \mathcal{L} \otimes \mathcal{J}^{\vee})$.  We would like to apply Bertini's theorem to these global sections to conclude that a general section gives a smooth divisor on $U^{sm}$.  By [Hart, III.10.9.2], the theorem (stated for projective varieties) holds for quasi-projective varieties when the system is finite-dimensional. Considering $X$ as an open subset of a projective variety $\overline{X}$, the inclusions $X \stackrel{i}{\rightarrow} \overline{X} \stackrel{j}{\rightarrow} \mathbb{P}$ determine $\mathcal{L} = i^*(j^*\mathcal{O}_{\mathbb{P}}(1))$ as the restriction of a very ample divisor on $\overline{X}$.  Similarly, the coherent sheaf $\mathcal{J}^{\vee}$ lifts to a coherent sheaf on $\overline{X}$, by [Hart, II.Ex.5.15].  Now on any projective variety, the global sections of a coherent sheaf form a finite dimensionsl vector space [Hart, II.5.19].  Putting this altogether, we may choose a very ample sheaf on $\overline{X}$ so that its tensor product with a lifting of $\mathcal{J}^{\vee}$ is globally generated.  The vector space of global sections is finite dimensional, hence so is the image of this vector space in $\Gamma(X, \mathcal{L} \otimes \mathcal{J}^{\vee})$.  The image generates the sheaf on $X$.  The same thing holds when restricting further to the open subscheme $U^{sm}$.  We denote by $V \subset \Gamma(X, \mathcal{L} \otimes \mathcal{J}^{\vee})$ this  finite-dimensional vector space of global sections that generate the sheaf $\mathcal{L} \otimes \mathcal{J}^{\vee}$.

Since $U^{sm}$ is smooth, Bertini's theorem implies that a general element of a finite-dimensional linear system without base points gives a smooth divisor on $U^{sm}$.  But the restriction of $V$ to $\Gamma(U^{sm}, (\mathcal{L} \otimes \mathcal{J}^{\vee})|_{U^{sm}})$ is such a system.  We conclude that there is a dense open subset of $\mathbb{P}(V)$ corresponding to divisors whose restrictions to $U^{sm}$ are smooth.

There are finitely many points of $X$ that are either generic points or singular codimension one points.  For each such point $p$, the subspace $W \subset V$ of sections vanishing at $p$ is a proper subspace because the sections in $V$ generate the sheaf.  The associated linear subspaces have as their union a proper subvariety of the projective space $\mathbb{P}(V)$, a finite union of vector subspaces is not the entire space, since the ground field is infinite).  Thus, we have two dense open subsets of $\mathbb{P}(V)$, which must therefore intersect.  Replacing $U$ by a smaller open subset (still containing the codimension one points) if necessary, it follows that there exist general global sections of $\mathcal{L} \otimes \mathcal{J}^{\vee}$ whose associated divisors are smooth and disjoint from the singular locus of $U$, after restriction to $U$.

Choosing such a section $s$, dualize $\mathcal{O}_X \stackrel{s}{\rightarrow} \mathcal{L} \otimes \mathcal{J}^{\vee}$.  Let $D$ be the cosupport of the resulting ideal sheaf $\mathcal{L}^{\vee} \otimes \mathcal{J}$.  Following Koll\'ar, there is a semiresolution $f: Y \rightarrow X$ such that the local models of $f^{-1}_*D + \textnormal{Ex}(f)$ are as in Section 2.  The inverse image ideal sheaf $(\mathcal{L}^{\vee} \otimes \mathcal{J}) \cdot \mathcal{O}_Y$ has cosupport equal to $f^{-1}D$.  Thus the inverse image ideal sheaf is invertible and has the correct normal form.

We do the same construction for the invertible ideal sheaf $\mathcal{L}^{\vee}$.  We can choose $\mathcal{L}$ so that it is generated by global sections.  Furthermore, we can choose a general section so that the associated divisor is smooth and disjoint from both Sing$(U)$ and the cosupport of $\mathcal{L}^{\vee} \otimes \mathcal{J}$.  This is possible once we have fixed a section of $\mathcal{L} \otimes \mathcal{J}^{\vee}$, since a general section of $\mathcal{L}$ is nonvanishing at a finite set of codimension one points.  Thus, after possibly shrinking $U$ (and keeping all the codimension one points), we can assume that $\mathcal{L}^{\vee}$ and $\mathcal{L}^{\vee} \otimes \mathcal{J}$ are ideal sheaves whose divisors are smooth and disjoint from one another and from Sing$(U)$ on $U$.  

Then the inverse images $$(\mathcal{L}^{\vee} \otimes \mathcal{J}) \cdot \mathcal{O}_Y = \mathcal{L}^{\vee} \cdot \mathcal{O}_Y \otimes \mathcal{J} \cdot \mathcal{O}_Y$$ and $\mathcal{L}^{\vee} \cdot \mathcal{O}_Y$ are both invertible.  Hence $\mathcal{J} \cdot \mathcal{O}_Y$ is invertible.  It corresponds to a (not necessarily effective) divisor $E$ whose components are smooth and intersect Ex$(f)$ and $C_Y$ transversally.\qed
\end{lemma}

We would like a semi log resolution of a pair $(X, I)$ when all of the summands of $I$ are divisorial sheaves.  If each of these is invertible in codimension one, then the same moving argument used in (4.1) implies that a semi log resolution exists.  In fact, when we include more than one divisorial sheaf in $I$, we can proceed inductively to obtain sections that do not vanish at any of the generic points, the singular codimension one points of $U$, or the codimension one points of inductively-defined divisors corresponding to components of $I$.  Note that for an arbitrary fractional ideal sheaf, we can always blow up to make the sheaf invertible, but the resulting morphism might not satisfy $f_*\mathcal{O}_Y = \mathcal{O}_X$.  

We say that the pair $(X, I)$ is \textit{effective} if the summands of $I$ are ideal sheaves and the coefficients are nonnegative rationals.

\begin{definition}  Let $(X, I)$ be an effective pair, and fix an integer $m \geq 2$.  Given a semi log resolution $f: Y \rightarrow X$ of $(X, I + \mathcal{O}_X(-mK_X))$, a boundary $\Delta$ on $X$ is said to be $m$-compatible for $(X, I)$ with respect to $f$ if: 

(i)  $m\Delta$ is integral and $\left\lfloor \Delta \right\rfloor = 0$,

(ii)  no component of $\Delta$ is contained in the cosupport of $I$,

(iii)  $f$ is a semi log resolution for $(X,\Delta; I + \mathcal{O}_X(-mK_X))$, and 

(iv)  $K^{\Delta}_{Y/X} = K_{m,Y/X}.$
\end{definition}

Once the sheaf $\mathcal{O}_X(-mK_X)$is invertible, pulling back commutes with composition of morphisms.  Our definitions of singularities need to be independent of the semiresolution chosen, and for that we need to be able to compose semiresolutions without affecting the pullback.

If there are $m$-compatible boundaries with respect to a semi log resolution of $(X, I + \mathcal{O}_X(-mK_X))$, then we say that the pair $(X, I)$ \textit{admits $m$-compatible boundaries}.  Before proving the existence of $m$-compatible boundaries, we prove the semismooth version of Bertini's theorem, which will be used in the proof that follows.

\begin{lemma}  Suppose $\mathcal{J}$ is an invertible sheaf on a semismooth variety $Y$, generated by a finite-dimensional subspace $V \subseteq \Gamma(Y, \mathcal{J})$.  Let $D$ be a given global normal crossing divisor on $Y$.  Then a general section in $V$ produces a divisor that forms global normal crossings with $D$.
\proof  Let $p: \overline{Y} \rightarrow Y$ be the normalization.  Since $\mathcal{J}$ is generated by global sections in $V$, there is a surjection $\bigoplus_{v \in V} \mathcal{O}_Y \rightarrow \mathcal{J}.$  Applying $p^*$ and following with the obvious surjection, we have the surjection $$\bigoplus_{p^*v} \mathcal{O}_{\overline{Y}} \rightarrow p^*\mathcal{J} \rightarrow \mathcal{J} \cdot \mathcal{O}_{\overline{Y}}.$$ Thus $\mathcal{J} \cdot \mathcal{O}_{\overline{Y}}$ is generated by the images of the $p^*v$. Given a general combination of vectors in $V$, we obtain a morphism $\mathcal{O}_Y \rightarrow \mathcal{J}$.  Then ``the same" combination of the pullbacks of those vectors produces the morphism obtained by applying $p^*$ and following with the obvious surjection.

Bertini's theorem for the smooth variety $\overline{Y}$ implies that a general combination of the $p^*v$ determines a divisor that makes simple normal crossings with the preimage of $D$.  Likewise, a general combination of the $v$ exhibits an effective divisor associated to $\mathcal{J}$.  Since the intersection of two general conditions is still general, we obtain a divisor $D'$ on $Y$ whose preimage makes simple normal crossings with the preimage of $D$.  (We need that the cosupport of $\mathcal{I} \cdot \mathcal{O}_{\overline{Y}}$ is just the preimage of the cosupport of $\mathcal{I}$, for any ideal sheaf $\mathcal{I}$.)  This implies that $D'$ makes global normal crossings with $D$.\qed
\end{lemma}

\begin{theorem}\textnormal{(Theorem 1.1)}  Every effective pair $(X,\mathcal{I})$ where the components of $\mathcal{I}$ are rank one reflexive sheaves and invertible in codimension one admits $m$-compatible boundaries for any $m \geq 2$.
\proof  We choose a very ample sheaf $\mathcal{L}$ such that $\mathcal{L}^{\vee} \otimes \omega_X \cong \mathcal{I}$ is an ideal sheaf.  By hypothesis, $\omega_X$ is invertible in codimension one, and hence the same is true for $\mathcal{I}$.  Let $D$ be the effective divisor associated to $\mathcal{I}$.  By (4.1), there is a semi log resolution of $(X, \mathcal{O}_X(-mK_X) + I + \mathcal{O}_X(-mD))$.  Then $\mathcal{O}_X(-mD) \cdot \mathcal{O}_Y$ is invertible.  Let $E = f^{^\natural}(mD)$ be the associated effective divisor.  Then $$f^{\natural}(mK_X) = f^{\natural}(mK_X -mD + mD) = f^{\natural}(mK_X -mD) + E,$$ and since $K_X -D$ is Cartier,  $$K_{m,Y/X} = K_Y - f^*(K_X - D) - \frac{1}{m}E.$$

Next, we choose a very ample sheaf $\mathcal{N}$ such that $\mathcal{N} \otimes \mathcal{O}_X(-mD)$ is globally generated.  A general section is also a section of $\mathcal{N}$ and produces an isomorphism $\mathcal{N}^{\vee} \otimes \mathcal{O}_X(mD) \cong \mathcal{J}$ for some ideal sheaf $\mathcal{J}$.  On the divisor level, we have $G = M + mD$, where $G$ is an effective Cartier divisor associated to $\mathcal{N}^{\vee}$ and $M$ is an effective divisor associated to $\mathcal{J}$.  Then $$f^*\mathcal{N}^{\vee} \cong \mathcal{N}^{\vee} \cdot \mathcal{O}_Y = (\mathcal{J} \cdot \mathcal{O}_Y) \otimes (\mathcal{O}_X(-mD) \cdot \mathcal{O}_Y).$$  Thus $\mathcal{J} \cdot \mathcal{O}_Y$ is invertible.  By definition of the inverse image of a divisorial sheaf, we have $(\mathcal{J} \cdot \mathcal{O}_Y)^{\vee} \cong \mathcal{J}^{\vee} \cdot \mathcal{O}_Y$.

It follows that $f^*G = f^{\natural}(M) + f^{\natural}(mD)$, where we have used the fact that additivity holds if both pullbacks are invertible.  Being an image of the sheaf $f^*(\mathcal{N} \otimes \mathcal{O}_X(-mD))$, $\mathcal{J}^{\vee} \cdot \mathcal{O}_Y$ is generated by global sections.  We conclude that $f^{\natural}(M)$ corresponds to a sheaf generated by its global sections.

Write $\mathcal{M} = \mathcal{J}^{\vee} \cdot \mathcal{O}_Y$, so that $\mathcal{M}$ is an invertible sheaf, generated by global sections.  As in the proof of (4.1), we may assume that $\mathcal{N}$ is the restriction of a very ample divisor on the projective closure of $X$, so that there is a finite-dimensional space of global sections generating $\mathcal{N} \otimes \mathcal{O}_X(-mD)$.  Take the inverse image $\mathcal{M} \cdot \mathcal{O}_{\overline{Y}}$ and apply Bertini's theorem.  A general section gives a divisor whose image in $Y$ makes global normal crossings with everything already in global normal crossing (including the conductor).  We can choose such a section generally enough that also the corresponding divisor in $\overline{Y}$ has no components in common with the preimage of $D$ and the preimages of the components of $I$.  Then in $X$ we necessarily have no components in common with $D$ or the components of $I$.

Therefore, we let $\Delta = \frac{1}{m}M$.  Then $m\Delta$ is integral and $\left\lfloor \Delta \right\rfloor = 0$ since $M$ can be chosen reduced.  We have $K_X + \Delta = K_X - D + \frac{1}{m}G$ is $\mathbb{Q}$-Cartier.  If again we choose a sufficiently general section and look at the normalization, then $f$ is a log resolution of the log pair $(X, \Delta; \mathcal{O}_X(-mK_X) + I)$.  The final condition for $m$-compatibility follows from the computation $$K^{\Delta}_{Y/X} = K_Y + \Delta_Y - f^*(K_X + \Delta)$$ $$= K_Y + \Delta_Y - f^*(K_X + \Delta - \frac{1}{m}G) - \frac{1}{m}f^*G$$ $$= K_Y - f^*(K_X-D) - \frac{1}{m}E = K_{m,Y/X}.$$  This completes the proof.\qed
\end{theorem}
\end{section}
\begin{section}{Singularities on Demi-Normal Varieties}
Let $f: Y \rightarrow X$ be a semiresolution of a demi-normal variety $X$ (possibly without $\mathbb{Q}$-Cartier canonical class), and let $E$ be an exceptional prime divisor.  For any $m \geq 0$, define $$a_{m,Z} := \textnormal{ord}_E(K_{m,Y/X}) + 1 - \textnormal{val}_E(Z).$$  The pair $(X, Z)$ is semi log canonical (resp., semi log terminal) if there exists and integer $m_0$ such that $a_{m_0, E} \geq 0$ (resp., $> 0$) for every exceptional prime divisor over $X$.  Define also $a_E(X,Z) : = \textnormal{ord}_E(K_{Y/X}) + 1 - \textnormal{val}_E(Z).$  The pair $(X,Z)$ is semi canonical if $a_E(X,Z) \geq 1$ for every exceptional prime divisor over $X$.  

The reason for the ``discrepancy" between the two definitions is that the important properties of the resulting class of singularities should be as in the case where $K_X$ is $\mathbb{Q}$-Cartier.  We have the following results, stated without proof.

A pair $(X,Z)$ admitting semi log resolutions is semi log canonical (resp., semi log terminal) if and only if there is a boundary $\Delta$ such that $(X,\Delta;Z)$ is semi log canonical (resp., semi log terminal).  In other words, we look at the order of $E$ in the log relative canonical divisor $K^{\Delta}_{Y/X}$ in place of $K_{m,Y/X}$.  

Suppose $f: Y \rightarrow X$ is a semiresolution of $X$, with the property that $\mathcal{O}_X(mK_X)\cdot \mathcal{O}_Y$ is invertible for some integer $m$.  Then $\textnormal{ord}_E(K^{-}_{Y/X}) > -1$ for every exceptional prime divisor $E$ of $f$ if and only if $X$ is semi log terminal.

A pair $(X,Z)$ is semi log canonical (resp., semi log terminal) only if the pair $(\overline{X}, Z \cdot \mathcal{O}_{\overline{X}})$ is log canonical (resp., log terminal) in the generalized sense of de Fernex and Hacon.

Let $X$ be demi-normal, and suppose that $Z = \Sigma a_k \cdot Z_k$ is an effective $\mathbb{Q}$-linear combination of effective Cartier divisors on $X$.  Then the pair $(X,Z)$ is semi canonical if and only if for all sufficiently divisible $m \geq 1$ (in particular, we ask that $ma_k \in \mathbb{Z}$ for every $k$), and for every semi log resolution of $(X, Z + \mathcal{O}_X(mK_X))$, there is an inclusion $$\mathcal{O}_X(m(K_X + Z)) \cdot \mathcal{O}_Y \subseteq \mathcal{O}_Y(m(K_Y + Z_Y))$$ as sub-$\mathcal{O}_Y$-modules of $\mathcal{K}_Y$, where $Z_Y$ is the birational transform of $Z$.

Suppose that $f: Y \rightarrow X$ is a semi log resolution of $(X,0)$.  If $\textnormal{ord}_E(K_{Y/X}) \geq 0$ for every exceptional divisor of $f$, then $X$ is semi canonical.  

If the log pair $(\overline{X}, C_{\overline{X}})$ is canonical, then $X$ is semi canonical.

The proofs of these facts are not difficult.  They either follow from the same proofs in \cite{dFH09}, or else from a comparison of $X$ with its normalization.  Note that when $f: Y \rightarrow X$ is a semiresolution, then the induced morphism on normalizations $\overline{Y} \rightarrow \overline{X}$ is a resolution of singularities.  

In the final section of this paper, we explore an example which shows that, when following the definitions given above, one can have a semi canonical demi-normal variety that is not semi log terminal.  This distinction evidently arises only in case the canonical divisor is not $\mathbb{Q}$-Cartier.  Whether or not it suggests a different set of definitions is necessary is not something we investigate.  We include it mainly because it illustrates the use of tools covered in the present paper.  Future research will be needed to show whether there is anything inherently wrong with the present definitions.
\end{section}
\begin{section}{An Example}
Suppose that $Y$ is a projectively-embedded variety whose canonical class is Cartier, and with semi canonical singularities.  The projectivized cone $X$ has $\mathbb{Q}$-Cartier canonical class if and only if $mK_Y = nH$ for some nonzero integers $m$ and $n$, where $H$ is the very ample divisor giving the embedding of $Y$.

Now let $X := \mathbb{P}^1 \times \mathbb{P}^1$ be embedded by the very ample sheaf $\mathcal{M} = \mathcal{O}_X(1,3)$.  We consider first a double covering of $X$ ramified over a general nonsingular curve, the curve being given by a general section of the sheaf $\mathcal{L}^{\otimes 2} $, where $\mathcal{L}= \mathcal{O}_X(0,2)$.  Such curves exist by Bertini's Theorem.  The double covering can be described as $$ W = \textbf{\textnormal{Spec}}_X(\mathcal{O}_X \oplus \mathcal{L}^{\vee}).$$  Next, we let $Y$ be a smooth curve in the linear system corresponding to $\mathcal{M}$, and we consider the induced double cover $$Z = \textbf{\textnormal{Spec}}_Y(\mathcal{O}_Y \oplus i^*\mathcal{L}^{\vee}).$$  There is a diagram $$\begin{CD}
Z         @>j>>      W \\
@VqVV                  @VpVV \\
Y         @>i>>      X \\
\end{CD}.$$     
Note that both $Z$ and $W$ are smooth, by \cite{KM98} 2.51, and that $j$ is a closed immersion because closed immersions are stable under base extension.  

The ample sheaf $p^*\mathcal{M}$ is in fact very ample on $W$.  By this result, we have a commutative diagram
$$\begin{CD}
W           @>p^*\mathcal{M}>>       \mathbb{P}^{M+N+1} \\
@VpVV                                               @VVV \\
X            @>\mathcal{M}>>             \mathbb{P}^M \\
\end{CD}.$$
All horizontal morphisms are closed immersions, all objects are smooth, and the hyperplane sections $Z$ and $Y$ are one-dimensional.  The rightmost vertical arrow is a linear projection.  Moreover, $W$ and $X$ embed as projectively normal varieties, since both very ample sheaves have degree six on a smooth variety.  In particular, the cone is normal at the vertex.

Since the diagram commutes, there is an induced diagram on the projectivized cones with respect to these embeddings, namely $$\begin{CD}
C(Z)    @>>>   C(W) \\
@VVV                @VVV \\
C(Y)   @>>>    C(X) \\
\end{CD}.$$
There is an involution $\tau$ on $Z$ that interchanges elements in the fibers of $Z$ over points of $Y$.  Since $Z$ and $Y$ are smooth projective curves and the base field has characteristic zero, $Z \rightarrow Y$ is ramified at only finitely many points, so the fixed point locus of $\tau$ is a finite set.  The fixed point locus is nonempty because $Y$ has intersection number 4 with the ramifying curve for $W \rightarrow X$.

There is an induced involution on $C(Z)$, obtained as follows.  Letting $\pi: C(Z) - P \rightarrow Z$ be the projection from the vertex $P$, we see that $Z$ is covered by open sets $U_i$ such that $\pi$ is given by $U_i  \times \mathbb{A}^1 \rightarrow U_i$.  In other words, the fibers of the projection are isomorphic to $\mathbb{A}^1$, and we send one line isomorphically to the other if the corresponding points of $Z$ are interchanged by $\tau$.  Sending the vertex to itself, we thus obtain an involution on $C(Z)$ whose fixed point locus is nonempty and of pure codimension one.

Now we can describe the variety whose singularities we want to compute.  It is the universal pushout $$\begin{CD}
C(Z)           @>>>   C(W) \\
@VVV                      @VVV \\
C(Z)/\tau  @>>>   V \\
\end{CD}.$$
We claim that $V$ has a canonical class which is not $\mathbb{Q}$-Cartier, and that $V$ is demi-normal.  Its normalization is $C(W)$.

We know that $C(W) \rightarrow V$ is finite and birational.  Moreover, $C(W)$ is normal because away from the vertex it is smooth (a resolution is given by blowing up the vertex and is isomorphic to a $\mathbb{P}^1$-bundle over $W$) and at the vertex it is normal because $W$ is projecitvely normal.  Thus $C(W)$ is the normalization of $V$.

That $V$ has a canonical class that is not $\mathbb{Q}$-Cartier follows from the fact that $C(W)$ has the same property.  This is due to the property we mentioned at the beginning of this section.  Note that the conductor is given by the hyperplane section $C(Z)$ and hence is a Cartier divisor.

To show that $V$ has the three properties $G_1$, $S_2$, and SN, we look at a similar pushout diagram involving $\mathbb{P}^1$-bundles over $Z$ and $W$:  $$\begin{CD}
\mathbb{P}(Z)         @>>>       \mathbb{P}(W) \\
@VVV                                         @VVV \\
\mathbb{P}(Z)/\tau @>>>       V' \\
\end{CD}.$$
The involution on $\mathbb{P}(Z)$ is obtained as before.  In particular, since $\mathbb{P}(Z)$ and $\mathbb{P}(W)$ are smooth and the fixed point locus has pure codimension one, $V'$ is semismooth.  The universal property of the pushout gives a morphism $f: V' \rightarrow V$, which is proper and birational.  With claim that $C(W)$ is $S_3$.  Then $C(Z)$ is $S_2$, hence normal.  By Zariski's Main Theorem, the proper birational morphisms $\mathbb{P}(Z) \rightarrow C(Z)$ (and likewise for $C(W)$) induce isomorphisms of structure sheaves.  Since the structure sheaf of the pushout is the pullback of the structure sheaves, we find that in fact $f_*\mathcal{O}_{V'} = \mathcal{O}_V.$  Since $V$ is $S_2$, this implies that $f$ is an isomorphism over codimension one points of $V$, and hence $V$ is also $G_1$ and SN (with $S_2$, SN in codimension one implies SN; see \cite{GT80}), because $V'$ is semismooth and such varieties are always demi-normal.

We now show that the normal variety $C(W)$ is not log terminal, and that the pair $(C(W),C(Z))$ is canonical.  It follows from the results presented in the previous section that $V$ is semi canonical but not semi log terminal.  

Recall that the canonical divisor on the projective bundle over over $W$ is given by $K_{\mathbb{P}} \sim \pi^*K_W -2W_0 + \pi^*(-L)$, where $L$ is the very ample divisor giving the embedding of $W$.  Pushing forward by the resolution $g: \mathbb{P}(W) \rightarrow C(W)$, we get $K_{C(W)} = f_*K_{\mathbb{P}} = C_{K_W} - C_L.$  On the cone, we are looking for a boundary $\Gamma = C_{\Delta}$ such that $K_{C(W)} + C_{\Delta}$ is $\mathbb{Q}$-Cartier.  However, a divisor on the cone is $\mathbb{Q}$-Cartier if and only if the associated divisor on $W$ is linearly equivalent to a rational multiple $kL$ of $L$.  This follows from a previous remark that a multiple of a hyperplane section of $W$ lifts to a multiple of a hyperplane section of $C(W)$.  Any other hypersurface section of $W$ will require two generators at the vertex.  Therefore, we are looking for effective divisors of the form $\Delta = sL - K_W$, where $s = k+1$.  

For semi log terminal singularities, we must compute the $m$-th limiting relative canonical divisor.  It suffices to calculate the order of $W_0$  in $K^-_{\mathbb{P}/C(W)}$.  We know that $K^{\Gamma}_{\mathbb{P}/C(W)} \leq K_{m,\mathbb{P}/C(W)}$ for all boundaries $\Gamma$ and all values of $m$, and by the main result, there exists a boundary divisor such that this is an equality, for any given $m$.  Thus we want to compute the supremum of the numbers $$\textnormal{ord}_{W_0}(K^{\Gamma}_{\mathbb{P}/C(W)}),$$ where $\Gamma$ is a boundary.  We calculate that $$K^{\Gamma}_{\mathbb{P}/C(W)} = K_{C(W)} + g^{-1}_*\Gamma - g^*(K_{C(W)} + \Gamma)$$ $$= \pi^*K_W - 2W_0 + \pi^*(-L) + g^{-1}_*\Gamma - \pi^*(K_W-L+\Delta) - (s-1)W_0$$ $$ = -(s+1)W_0 + g^{-1}_*\Gamma - \pi^*\Delta.$$  So, if $t$ is the infimum of values $s$ such that $sL-K_W$ is effective, then the value of the relative canonical divisor at $W_0$ is $-(1+t).$

For semi canonical singularities, we need to calculate the valuation along $W_0$ of the relative canonical divisor $K_{\mathbb{P}/C(W)}$.  The proof of the main theorem can be modified to show that $$\textnormal{val}_{W_0}(K_{\mathbb{P}/C(W)}) = \textnormal{sup}\{\textnormal{ord}_{W_0}(K_{\mathbb{P}} + g^*(-K_{C(W)}+ \Gamma'))\},$$ where $-K_{C(W)} + \Gamma'$ is $\mathbb{Q}$-Cartier and $\Gamma'$ is effective.  Here $\Gamma'$ corresponds to $\Delta'=rL+K_W$ with $r$ rational, which is the condition for divisors on the cone to be $\mathbb{Q}$-Cartier.  Thus if $t$ is the smallest value of $r$ such that $rL+K_W$ is effective, then the valuation of $K_{\mathbb{P}/C(W)}$ along $W_0$ will be $t-1$.

Let $L$ be the divisor corresponding to the sheaf $\mathcal{M}$.  To determine whether $C(W)$ has log terminal singularities, we need to compute the smallest value of $s$ such that $sp^*L-K_W$ is linearly equivalent to an effective divisor.  If $s=0$, we have $-K_W$, which corresponds to the sheaf $p^*\mathcal{O}_X(2,0)$, which clearly has global sections.  In fact, its pushforward is $\mathcal{O}_X(2,0) \oplus \mathcal{O}_X(2,-2)$.  Any nonzero global section of the form $(t,0)$ will work.  For negative $s$, $sp^*L-K_W$ is linearly equivalent to an effective divisor if and only if one of its positive integer multiples is.  The corresponding sheaf is $p^*\mathcal{O}_X(r(s+2), 3rs)$ for a positive integer $r$.  The pushforward is $\mathcal{O}_X(r(s+2), 3rs) \oplus \mathcal{O}_X(r(s+2), 3rs-2)$.  Neither summand has nonzero global sections, since the second factors are negative.  The valuation of the relative canonical divisor $K^-_{\mathbb{P}/C(W)}$ along the exceptional divisor is therefore $-(1+s) = -1$.  In particular, $C(W)$ is not log terminal.

For canonical singularities, we look for the smallest value of $r$ such that $rp^*L+K_W$ is effective.  The corresponding sheaf is $\mathcal{O}_X(r-2,3r)$.  As in the previous part of the proof, we see that when $r=2$, there are nonzero global sections, while if $r<2$ there can be none.  Thus the relative canonical divisor has value $r-1=1$ along the exceptional divisor.  We claim that the value of the conductor $C(Z)$ along the exceptional divisor is 1.  In fact, $C(Z)$ is a hyperplane section of $C(W)$, so its multiplicity at the vertex is exactly one.  When we pull back via $g$, we get exactly one copy of the exceptional divisor.  Then it follows that the pair $(C(W),C(Z))$ is canonical.

\end{section}
\nocite{Art70}
\nocite{Kol13}
\nocite{Hart94}
\nocite{Hart87}
\nocite{BH98}
\nocite{GT80}
\nocite{Kov99}
\nocite{LV81}
\nocite{Reid94}
\nocite{dFH09}
\nocite{KSS09}
\nocite{Sch09}
\nocite{vS87}
\nocite{Berq14}
\nocite{GR70}
\nocite{Berq14b}
\nocite{KM98}
\nocite{K00}

\bibliographystyle{plain}
\bibliography{paper}

\end{document}